\documentclass[10pt]{article}
\usepackage{amsmath}
\usepackage{amsfonts}
\usepackage{graphicx}
\usepackage{latexsym, amscd, euscript}
\usepackage{color}
\usepackage{amssymb}
\usepackage{makeidx}
\usepackage{amsthm}
\usepackage{epsfig}

\author{Fuensanta Aroca\thanks{Research partially supported by PAPIIT-UNAM, CONACYT (Mexico) grants 55084,
162340 and 117110}\\(Universidad Nacional Aut\'onoma de M\'exico.)
\\ Mirna G\'omez-Morales.\thanks{Research supported CONACYT (Mexico) grant 213186‏}
(University of Warwick.)
\\ Khurram Shabbir\thanks{Partially supported by CONACYT (Mexico) and TWAS (Italy) grant FR
3240223595. The third author wants to thank these institutions for their support.}\\(Department of Mathematics, GC University Lahore.)}

\title{Torical Modification of Newton non-degenerate ideals.}

\begin{document}
\hyphenation{para-me-trizations}
\newtheorem{teo}{Theorem}[section]
\newtheorem{coro}[teo]{Corolary}
\newtheorem{lema}[teo]{Lemma}
\newtheorem{prop}[teo]{Proposition}
\newtheorem{obs}[teo]{Observation}
\newtheorem{af}{Afirmation}
\theoremstyle{definition}
\newtheorem{defi}[teo]{Definition}
\newtheorem{Remark}[teo]{Remark}
\newcommand{\C}{{\mathbb C}}
\newcommand{\K}{{\mathbb K}}
\newcommand{\N}{{\mathbb N}}
\newcommand{\Q}{{\mathbb Q}}\printindex
\newcommand{\R}{{\mathbb R}}
\newcommand{\Z}{{\mathbb Z}}
\newcommand{\sA}{{\mathcal A}}
\newcommand{\vertice}{{\bf P}}
\maketitle

 {\em Dedicated to professor Hironaka in his 80th birthday.}

\section{Introduction.}

Newton non degenerate singularities were introduced in the 70's\cite{Ehlers:1975,Kouchirenko:1976,Khovanskii:1977b}.
For hypersurfaces, the definition is given in terms of the Newton Polyhedron of the function defining the hypersurface.
For complete intersection singularities they are characterized in terms of the Newton polyhedra of a given set of generators of the ideal.
The reference book on the subject is \cite{Oka:1997}.

Newton non degenerate singularities of hypersurfaces and of complete intersections have been widely studied (See for example \cite{BraunNemethi:2007,Kouchirenko:1976} and  \cite{Khovanskii:1977,Oka:1990,Bivia:2007,SaiaZuniga:2005}). A good resolution of a non degenerate singularity may be constructed from the dual fan of the Newton boundaries.

In this paper we extend the definition of Newton non degenerate to non necessarily complete intersection singularities. Our definition is new and does not depend on the system of generators chosen and is given in terms of initial ideals.

The Gr\"{o}bner fan of an ideal is the extension to non-principal ideals of the concept of fan dual to the Newton polyhedron. The tropical variety associated to a hypersurface $H$ is the ($dim\, H-1$)-skeleton of the fan dual to the Newton polyhedron of the function defining the hypersurface. Reference books on Tropical Geometry are \cite{Gathmann:2006,ItenbergMikhalkin:2007,RichterSturmfels:2005}.

In \cite{ArocaIlardiLopez:2010} an extension of the Newton-Puiseux method to compute parametrizations of plane curves is extended to arbitrary codimension replacing the Newton Polyhedron by the tropical variety.

We prove that a regular refinement of the Gr\"{o}bner fan of the ideal defining the non degenerate variety gives a resolution. We also prove that the strict transform intersects the orbit associated to a cone if and only if the cone is contained in the tropical variety. Both results are original statements.

The second author would like to thank Carles Bibia-Ausina and Dmitry Kerner, for clarifying discussions on the subject and, to Meral Tosun, for the example in section \ref{Variedades no degeneradas}.
 The third author would like to thank A. Jensen for answering several questions by e-mail. In particular for the proof of Proposition \ref{Prop5.161}.

\section{Cones and fans.}

In this section we introduce some basic concepts of convex geometry. These concepts may be found in several books (see for example \cite{Fulton:1993}).

Given vectors $u^{(1)}, \ldots, u^{(k)} \in \R^n$. The \textbf{polyhedral cone}  generated by $u^{(1)}, \ldots, u^{(k)}$ is the set
\[
\sigma =\langle u^{(1)}, \ldots, u^{(k)} \rangle :=\{\lambda_1 u^{(1)}+\cdots+ \lambda_k u^{(k)}; \;\lambda_i\in \R_{\geq 0},\; i=1, \ldots, k\}\subset \R^n.
\]
The vectors $u^{(i)}$'s are called the \textbf{generators} of the cone. A polyhedral cone is said to be  \textbf{rational} if it has a set of generators in $\Z^n$.

We will denote by $e^{(1)},\ldots , e^{(n)}$ the vectors in the standard base of $\R^n$. With this notation the first orthant is the cone $(\R_{\geq 0})^n =\langle e^{(1)},\ldots , e^{(n)}\rangle$. 

The cone generated by the columns of a matrix $M$ will be denoted by $\langle M \rangle$.

The \textbf{dimension} of $\sigma$ is the dimension of the minimal linear subspace ${\mathcal L} (\sigma)$ containing $\sigma$ and is denoted by $dim(\sigma)$. The dimension of $\langle M \rangle$ is equal to the rank of the matrix $M$.

The \textbf{relative interior} of a cone $\sigma$ is the interior of $\sigma$ as a subset of ${\mathcal L} (\sigma )$. That is
$$Int_{rel} \langle u^{(1)},\ldots, u^{(s)}\rangle = \{ \lambda_1 u^{(1)} +\cdots + \lambda _s u^{(s)}; \, \lambda_i \in \R_{>0}\}.$$

The \textbf {dual} $\sigma^{\surd}$ of a cone $\sigma$ is the cone given by 
$$\sigma^{\surd} :=  \{ v \in \R^n ; v \cdot u \geq 0,\, \text{for all}\, u\in \sigma \}$$
where $u\cdot v$ stands for the inner product $u\cdot v := u_1v_1+\cdots +u_nv_n$  of the vectors$u=(u_1,\ldots ,u_n)$ and $v=(v_1,\ldots ,v_n)$.

Let $M$ be an unimodular matrix. We have
\begin{equation} \label{generadores del dual de regular}
{\langle M \rangle}^{\surd}= \langle {\left( M^{-1}\right) }^t \rangle
\end{equation}
where $M^t$ stands for the transpose of the matrix $M$.

A rational polyhedral cone is said to be \textbf{strongly convex} if it does not contain any non-trivial linear subspace. A cone contained in the first orthant is strongly convex. The dual of a cone of maximal dimension is strongly convex.

A vector in $\Z^n$ is said to be \textbf{primitive} when the maximum common divisor of its coordinates is 1. 

The set of vectors $\{u^{(1)}, \ldots, u^{(k)} \} \subset \Z^n$ is the set of \textbf{vertices} of a rational strongly convex cone $\sigma$ when
\begin{itemize}
	\item
	$u^{(i)}$ is primitive for $i\in \{ 1,\ldots ,k\}$.
	\item
	$\sigma = \langle u^{(1)}, \ldots, u^{(k)}\rangle$.
	\item
	$\langle u^{(1)},\ldots ,u^{(i-1)},u^{(i+1)},\ldots ,u^{(k)}\rangle\subsetneq \sigma$ for $i=1,\ldots ,k$.
\end{itemize}

By $\sigma=Cone(u^{(1)}, \ldots, u^{(k)} )$ we will denote the rational convex cone 
 with vertices $u^{(1)}, \ldots, u^{(k)}$. We will also write $\sigma = Cone (M)$ where $M$ is the matrix that has as columns the vertices of $\sigma$.

A strongly convex rational polyhedral cone $\sigma=Cone ( u^{(1)}, \ldots, u^{(k)})$ is said to be \textbf{regular} when the group ${\mathcal L} (\sigma )\cap \Z^n$ is of rank $k$ and is generated by the vertices of $\sigma$.

\begin{Remark}
The vertices of a regular cone are linearly independent over $\Q$. An $n$-dimensional rational cone in $\R^n$ is regular if and only if $\sigma =Cone (M)$ where $M\in GL (\Z ,n)$ is an unimodular matrix.
\end{Remark}

\begin{Remark}\label{caras de cono regular}
 Let $\sigma =Cone (u^{(1)}, \ldots, u^{(k)} )$ be a regular cone. The faces of $\sigma$ are the cones $Cone (u^{(i_1)}, \ldots, u^{(i_s)} )$ with $\{ i_1,\ldots ,i_s\}\subset\{1,\ldots ,k\}$.
\end{Remark}

A collection of cones $\Sigma$ in $\R^n$ is called a \textbf{polyhedral fan} if it satisfies the following properties:
    \begin{itemize}
        \item[i)] Every face of a cone in $\Sigma$ is a cone in $\Sigma$;
        \item[ii)] The intersection of any two cones $\sigma, \tau \in \Sigma$ is a face of both $\sigma$ and $\tau$.
    \end{itemize}

 A polyhedral fan is said to be \textbf{regular} if all of its cones are regular. 

\begin{Remark}\label{Matriz asociada a una cara de cono regular}
	Let $\Sigma$ be a regular fan and let $\tau$ and $\sigma$ be cones in $\Sigma$. By remark \ref{caras de cono regular}, the cone $\tau$ is a face of $\sigma$ if and only if the set of vertices of $\tau$ is a subset of the set of vertices of $\sigma$. That is $\tau =Cone (T)$ and $\sigma = Cone (M)$ where $M= (T,S)$.
\end{Remark}

A fan $\Sigma^{'}$ is a \textbf{refinement} of a fan $\Sigma$ if every $ \sigma \in \Sigma$ is a union of cones in $\Sigma^{'}$. A refinement $\Sigma^{'}$ is called \textbf{regular} if every cone in $\Sigma$ is regular.

\begin{prop}
 Any fan has a regular refinement.
\end{prop}
The proof is left as an exercise in section 2.6 of \cite{Fulton:1993}.

\section{Newton Polyhedron.}

Let $\K$ be an algebraically closed field of any characteristic.

A  polynomial in $n$ variables with coefficients in $\K$ is an expression of the form
\begin{eqnarray} 
\label{Pol}
f(x)=\sum_{\mu\in \Omega \subset  \Z_{\geq 0}^n} a_{\mu}x^{\mu};\,\, \# \Omega <\infty, a_{\mu} \in \K
 \end{eqnarray}
 where $x^{\mu}:=x_1^{\mu_1}\cdots x_n^{\mu_n}$.

The ring of polynomials in $n$ variables with coefficients in $\K$ will be denoted by $\K [x_1,\ldots ,x_n]$.

The \textbf{support} or \textbf{set of exponents} of $f$ is defined by \begin{eqnarray*}
\varepsilon (f):=\{\mu\in {\Z_{\geq 0}}^n;  a_{\mu}\neq 0\}.
 \end{eqnarray*}


\begin{Remark}\label{Si f no se anula el origen es un exponente}
	Given $f \in \K [x_1,\ldots ,x_n]$, as in (\ref{Pol}),
 we have $f(\underline{0})= a_{\underline{0}}$, and, then,
	 \[
	f (\underline{0}) \neq 0 \Longleftrightarrow \underline{0} \in \varepsilon (f).
	\]
\end{Remark}

The \textbf{Newton polyhedron} of $f \in \K [x_1,\ldots ,x_n]$ is the convex hull
\begin{eqnarray*}
	NP(f):=Conv(\{\mu + (\R_{\geq 0})^n ; \mu \in \varepsilon (f)\} ).
\end{eqnarray*}

\begin{figure}[h]
\begin{center}
\includegraphics[totalheight=1.7cm]{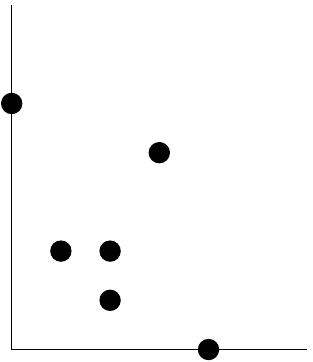}\,\,\,\,\,\,\,\,\,\,\,\,\,\,\,\,\,\,\,\,\,\,\,\,\includegraphics[totalheight=1.7cm]{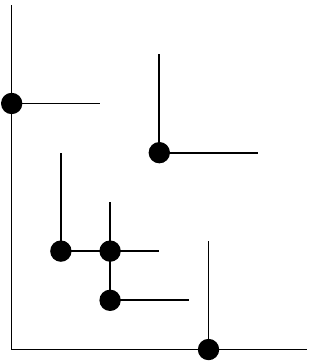}\,\,\,\,\,\,\,\,\,\,\,\,\,\,\,\,\,\,\,\,\,\,\,\, \includegraphics[totalheight=1.7cm]{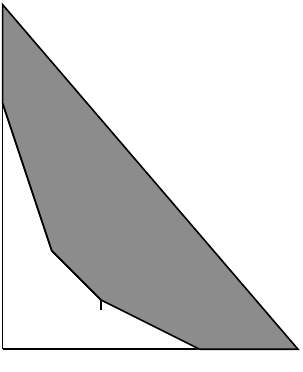}
\caption{$\varepsilon(f),\varepsilon(f)+\R^2_{\geq0},NP(f)$}

\end{center}
\end{figure}

\begin{Remark}\label{NP of a unit}
Given $f \in \K[x_1,\ldots ,x_n]$, by Remark \ref{Si f no se anula el origen es un exponente}, $f (\underline{0}) \neq 0$ if and only if the Newton polyhedron of $f$ is the first orthant.
\end{Remark}

\begin{Remark}\label{Cuando el poliedro tiene solo unvertice es un monomio veces una unidad}
Let $f$ be a polynomial in $\K [x_1,\ldots ,x_n]$. The Newton polyhedron of $f$  has only one vertex if and only if
\[
f=x^\alpha h
\] 
where $h$ is a polynomial in $\K [x_1,\ldots ,x_n]$ with $h(\underline{0})\neq 0$, and $\alpha$ are the coordinates of the vertex of $NP (f)$.
\end{Remark}

\begin{figure}[h]
\begin{center}
\includegraphics[totalheight=4.7cm]{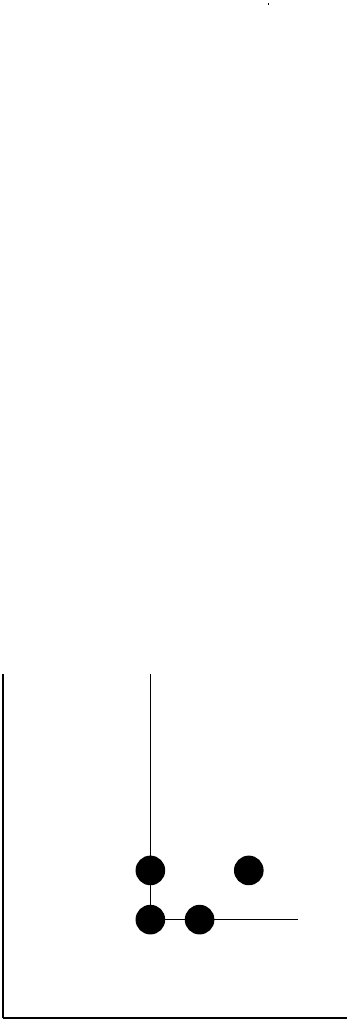}\,\,\,\,\,\,\,\,\,\,\,\,\,\,\,\,\,\,\,\,\,\,\,\,\,\,\,\,\,\,\,\,\,\,\,\,\,\,\,\,\,\,\,\,\,\,\,\,\,\, \includegraphics[totalheight=4.7cm]{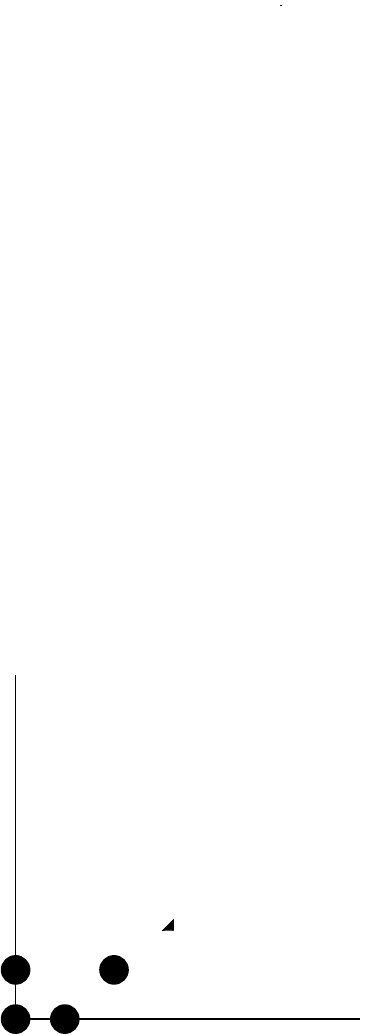}
\caption{$NP(f)$,\hspace{3cm} $NP(h)$ where $f=x^3y^2 h$}

\end{center}
\end{figure}

Let $f$ be a polynomial as in (\ref{Pol}) and let $F$ be a face of the Newton polyhedron of $f$. The {\bf restriction} of $f$ to the set $F\subset\Z^n$ is defined as
\[
f|_F := \sum_{\mu\in \varepsilon (f)\cap F \subset \Z^n} a_{\mu}x^{\mu}.
\]
\begin{Remark}\label{restingir es como poner ceros}
Let $f$ be a polynomial in $\K [x_1,\ldots ,x_n]$. We have
\[
f|_{\langle e^{(1)},\ldots ,e^{(s)}\rangle} = f(x_1,\ldots ,x_s,0,\ldots ,0).
\]
\end{Remark}

\section{The dual fan.}

Given $\omega\in (\R_{\geq 0})^n$ the $\omega$\textbf{-order} of a non-zero polynomial in $\K [x_1,\ldots ,x_n]$ is defined as
\begin{eqnarray}\label{omegaorden}
	\nu_{\omega}(f):=min\{\omega \cdot \mu; \mu\in \varepsilon(f)\}.
\end{eqnarray}

Let $f\in\K [x_1, \ldots , x_n]$ be a polynomial. Given a vector $\omega$ in the first orthant set
\[
 \pi_\omega (f): = \{ x\in\R^n ; \omega\cdot x =\nu_\omega (f)\}.
\]
The hyperplane $\pi_\omega (f)$ is a supporting hyperplane for $NP (f)$. 
\begin{figure}[h]
\centering
\includegraphics[totalheight=2.7cm]{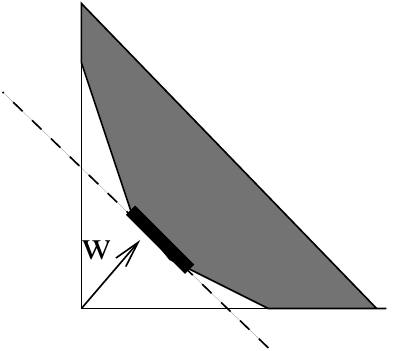}\,\,\,\,\,\,\,\,\,\,\,\,\,\,\,\,\,\,\,\,\,\,\,\,\includegraphics[totalheight=2.7cm]{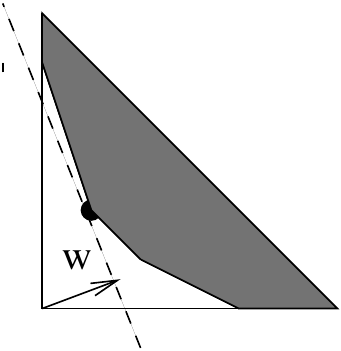}
\caption{$face_{w}(f)$}
\end{figure}
\medskip

The intersection
\[
face_{\omega} (f):= \pi_\omega (f)\cap NP (f)
\]
is a face of $NP (f)$. 

Given a face $F$ of $NP (f)$ set
\[
C_F := \{\omega \in  (\R_{\geq 0})^n ; F\subset face_\omega (f)\}.
\]

\begin{figure}[h]
\centering
\includegraphics[totalheight=1.7cm]{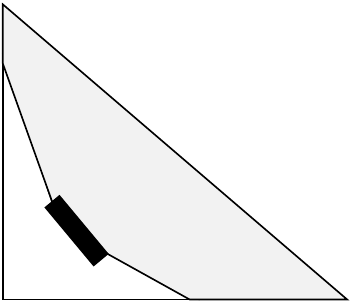}\,\,\,\,\,\,\,\,\,\,\,\,\,\,\,\,\,\,\,\,\,\,\,\, \includegraphics[totalheight=1.7cm]{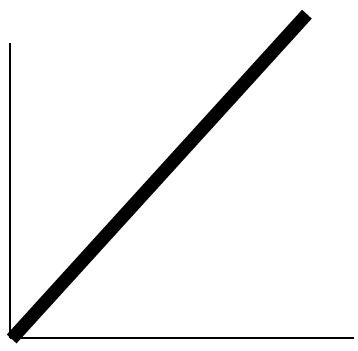}
\caption{$F,C_F$}
\includegraphics[totalheight=1.7cm]{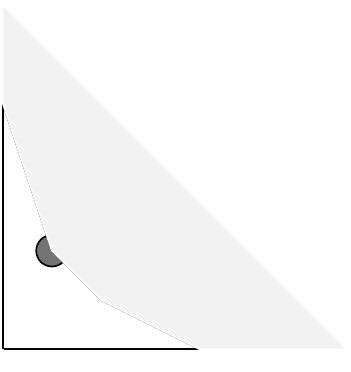}\,\,\,\,\,\,\,\,\,\,\,\,\,\,\,\,\,\,\,\,\,\,\,\,\includegraphics[totalheight=1.7cm]{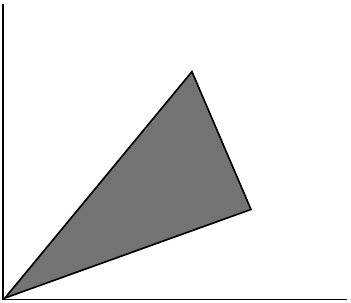}
\caption{$\vertice,C_\vertice$}

\end{figure}
\medskip

The collection of cones
\[
\Sigma(f):=  \{ C_F; F \,\text{is a face of}\,\, NP (f)\}
\]
forms a fan.

\begin{figure}[h]
\centering
\includegraphics[totalheight=1.7cm]{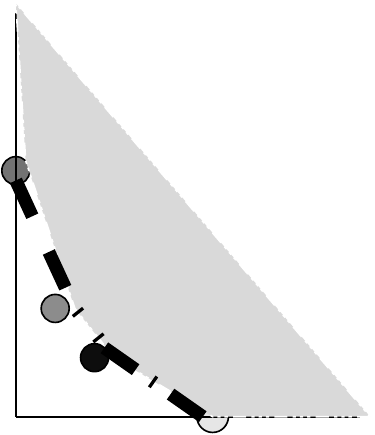}\,\,\,\,\,\,\,\,\,\,\,\,\,\,\,\,\,\,\,\,\,\,\,\,\,\,\,  \includegraphics[totalheight=1.7cm]{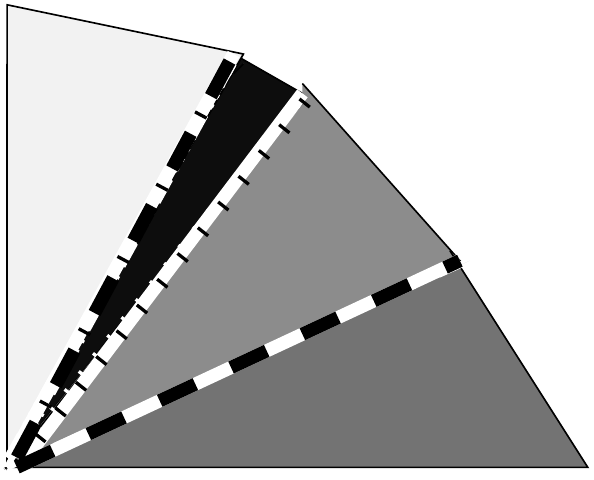}
\caption{}
\end{figure}
\medskip

\begin{Remark}\label{abanico de una unidad}
Given $f \in \K[x_1,\ldots ,x_n]$, if $f (\underline{0}) \neq 0$ then, by remark \ref{NP of a unit}, we have
\[
\Sigma (f)=\{\langle e^{(i_1)},\ldots , e^{(i_s)}\rangle ; \{i_1,\ldots ,i_s\}\subset \{1,\ldots ,n\}\}.
\]
\end{Remark}4

The mapping
\[
\begin{array}{ccc}
	\{\text{faces of } NP (f)\}	
		& \longrightarrow
			& \Sigma (f)\\
	F	& \mapsto 
			& C_F
\end{array}
\]
gives a duality.

For $\omega$ in the relative interior of the cone $C_F$ we have $face_\omega (f)= F$.

\begin{Remark}\label{relacion de equiv en principales}
	Given $\omega$  and $\omega'$ in the relative interior of the cone $\langle e^{(i_1)},\ldots ,e^{(i_s)}\rangle$.  The equality $face_{\omega}f= face_{\omega^{'}}f$ holds if and only if $\omega$ and $\omega'$ belong to the relative interior of the same cone of $\Sigma (f)$.
\end{Remark}

\begin{Remark} 
	Let $\vertice$ be a vertex of $NP (f)$. We have
	\[
	\varepsilon (f)\subset \vertice+ {C_\vertice}^\surd.
	\]
\end{Remark}

\begin{figure}[h]
\centering
\includegraphics[totalheight=1.7cm]{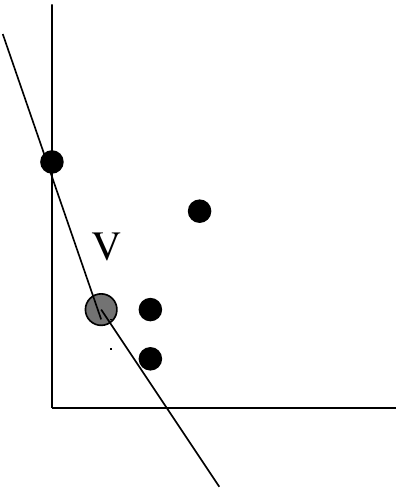} \,\,\,\,\,\,\,\,\,\,\,\,\,\,\,\,\,\,\,\,\,\,\,\,\,\,\, \includegraphics[totalheight=1.7cm]{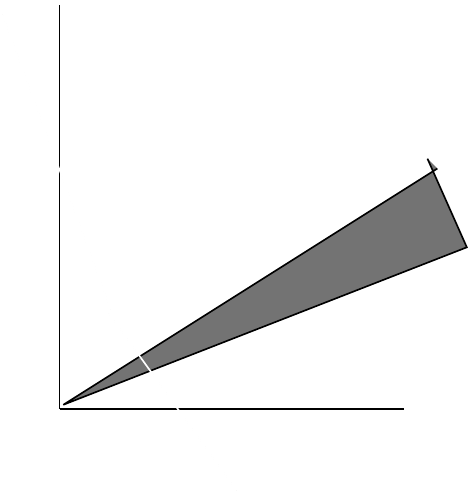}
\caption{$\vertice,C_\vertice$}
\end{figure}
\medskip

\begin{defi}\label{defi3.3.}
 	We will say that $\sigma$ is a \textbf{good cone} for $f$ if it is contained in a cone of $\Sigma (f)$.
\end{defi}

Let $f \in \K[x_1, \ldots ,x_n]$ be a polynomial and let $\sigma$ be a good cone for $f$. We have that $face_\omega f= face_{\omega '}f$ for $\omega,\omega'\in Int_{rel}\sigma$. We define the \textbf{$\sigma$-face} of $f$ as 
\[
face_{\sigma} f := face_{\omega} f\,\,\,\,\,\text{where}\,\,\,\omega\in Int_{rel}(\sigma).
\]

\begin{Remark} \label{El conjunto de exponentes esta contenido en el dual de un buen cono}
	Let $\sigma$ be a good cone for $f$ of maximal dimension. We have
	\[
	\varepsilon (f)\subset \vertice+ {\sigma}^\surd
	\]
where $\vertice$ is the vertex dual to the cone of $\Sigma (f)$ containing $\sigma$.
\end{Remark}

\section{Monomial transformations}

Given an unimodular matrix with integer entries,  $M\in GL(n,\Z)$, we will denote by $\phi_M$ the  morphism given by
\begin{center}
$\begin{array}{llll}
\phi_M :& (\K^{*})^n & \longrightarrow & (\K^{*})^n\\
 & z& \longmapsto & ({z}^{u^{(1)}}, {z}^{u^{(2)}},\ldots , {z}^{u^{(n)}})
\end{array} $
\end{center}
where $u^{(1)},u^{(2)},\ldots,u^{(n)} $ are the columns of the matrix $M$ and $\K^*:=\K\setminus\{ 0\}$.

The morphism $\phi_M$ is bi-rational on $\K^n$. It is bi-regular on the torus $(\K^{*})^n \simeq \mathbb{T} ^n$ and the following equalities are satisfied: \\
\begin{eqnarray}
 \phi_M \circ \phi_{M'} = \phi_{{M'}M } \\ 
 (\phi_M )^{-1} = \phi_{M^{-1}}.
\end{eqnarray}

Given a vector $\mu\in\Z^n$, we have
\begin{equation}\label{transformacion de un monomio}
	{\phi_M (x)}^\mu = 
	x^{\mu_1u^{(1)}}\cdots  x^{\mu_nu^{(n)}} = {x_1}^{\sum_{i=1}^n\mu_i u^{(i)}_1}\cdots {x_n}^{\sum_{i=1}^n\mu_i u^{(i)}_n}= x^{M\mu}.
\end{equation}

  We will denote by $L_M,$ the linear map given by
\[
\begin{array}{llll}
	L_{M} : & \R^n & \longrightarrow & \R^n\\
	 & x & \longmapsto & M x.
	\end{array}
\]

Let $f$ be a polynomial in $n$ variables as in (\ref{Pol}), equation (\ref{transformacion de un monomio}) implies
\begin{eqnarray}\label{fcomppi}
        f\circ \phi_M (x)= \sum_{\mu\in \varepsilon(f)} a_{\mu}x^{M\mu}
    \end{eqnarray}
then 
\begin{eqnarray}  \label{La transformacion en los exponentes}
  \varepsilon(f\circ\phi_M)= M \cdot \varepsilon(f)=L_M (\varepsilon(f)). 
\end{eqnarray}

Let $\pi$ be a supporting hyperplane for $NP (f)$, then $L_M (\pi)$ is a supporting hyperplane for $NP (f\circ \phi_M)$, and 
\begin{equation}\label{restriccion a hiperplanos}
{f\circ \phi_M}|_{L_M (\pi)\cap NP (f\circ \phi_M)} = f|_{\pi\cap NP (f)}\circ \phi_M.
\end{equation}

\begin{lema} \label{L sub M del dual es el primer ortante}
Let $M$ be an unimodular matrix with integer entries, we have
\[
L_M({(cone(M^t))}^{\surd})= (\R_{\geq 0})^n.
\]
 \end{lema}
Proof. 

The dual 
${cone(M^t)}^{\surd}$ is generated by $\langle v^{(1)}, \ldots, v^{(n)} \rangle,$ where by (\ref{generadores del dual de regular}), $v^{(1)}, \ldots, v^{(n)}$ are the columns of the matrix $M^{-1}$.

\[
\begin{array}{cl} 
	L_{M}({(cone(M^t))}^{\surd})
		&= \{L_{M} (t_1v^{(1)}+\cdots+ t_nv^{(n)})\,;\, t_i \in \R_{\geq0}\}\\
		&\{t_1 (M v^{(1)})+\cdots+t_n( M v^{(n)})\,;\, t_i \in \R_{\geq0}\}\\
		&=\{t_1(1,0,\ldots,0)+\cdots+t_n (0,\ldots,0,1)\,;\, t_i \in \R_{\geq0}\}\\
		& =\{(t_1,\ldots,t_n)\,;\, t_i \in \R_{\geq0}\}= (\R_{\geq 0})^n.
\end{array}
\]
Hence we have the required result.

\begin{prop} \label{El poliedro se convierte en el primer cuadrante trasladado}
  Let $M\in GL(n,\Z) $ be such that $\sigma= cone(M^t)$ is good for $f$, then  $NP (f\circ \phi_M)$  has only one vertex (i. e.  ${\R}^n_{\geq 0}$ is good for $f\circ \phi_M$).
\end{prop}
Proof. 

Suppose that $\sigma =cone (M^t)$ is good for $f$, then, by Remark \ref{El conjunto de exponentes esta contenido en el dual de un buen cono}, 
\begin{equation}\label{primera ecuacion en prueba de prop 5.2}
	\varepsilon (f)\subset \vertice+ \sigma^\surd= \vertice + (cone (M^t))^\surd
\end{equation}
with $\vertice\in\varepsilon (f)$.

By (\ref{La transformacion en los exponentes}) and (\ref{primera ecuacion en prueba de prop 5.2}),
\[
\varepsilon (f\circ \phi_M)\subset L_M(\vertice) + L_M((cone (M^t))^\surd)     =    L_M(\vertice) + {(\R_{\geq 0})}^n
\]
where the equality follows from Lemma  \ref{L sub M del dual es el primer ortante}.

Now, since $\vertice$ is in the support of $f$, the point  $L_M(\vertice)$ is in the support of  $f\circ \phi_M$.
Then, $L_M (\vertice)$ is the only vertex of $NP (f\circ \phi_M )$, and we have the result.

\begin{prop}\label{Como cambian las caras a las que hay que restringir}
	Let $M\in GL(n,\Z) $ be such that $\sigma= Cone(M)$ is good for $f$. Given $\lambda\subset \{1,\ldots ,n\}$, let $\tau$ be the face of $\sigma$ generated by the $i^{th}$ columns of $M$, with $i\in\lambda$. We have
	\[
	f|_{face_{\tau} (f)}\circ \phi_{M^t} = f\circ\phi_{M^t}|_{\vertice+ Cone (\{ e^{(i)}\}_{i\in\lambda^C})}
	\]
where $\vertice$ is the only vertex of $NP (f\circ\phi_{M^t})$.
\end{prop}

Proof. 

Consequence of equation (\ref{restriccion a hiperplanos}) and Proposition \ref{El poliedro se convierte en el primer cuadrante trasladado}.

\section{Toric modification.}

In this section we will recall how to construct the modification associated to a regular fan $\Sigma$ with support in the first orthant. See for example \cite[Cap.2]{Oka:1997}.

Let $\mathfrak{M}$ be the set of cones in $\Sigma$ of maximal dimension. Let $\sigma\in \Sigma$ be a cone in $\mathfrak{M}$.
 We will associate to $\sigma$ one copy of the affine space $\K^n$ and we will denote it by $U_{\sigma}$. Consider the disjoint union
\begin{eqnarray}
\mathfrak{C}= \bigsqcup_{{\sigma} \in \mathfrak{M}} U_{\sigma}.
\end{eqnarray}

Let $\sigma=Cone(M)$ and $\sigma'=Cone(M')$ be cones in $\mathfrak{M}$.

Consider 
\[
\phi_{(M^{-1})^t} : (\K^{*})^n \longrightarrow U_{\sigma}\quad \text{ and }\quad  \phi_{(M'^{-1})^t} : (\K^{*})^n \longrightarrow U_{\sigma'}.
\]

The composition 
$ \phi_{(M'^{-1})^t} \circ \phi^{-1}_{(M^{-1})^t} = \phi_{M^t.(M'^{-1})^t} :U_{\sigma}\dashrightarrow  U_{\sigma'} $  is a bi-rational morphism.

An equivalence relation is defined in $\mathfrak{C}$ as follows: Given two points $\;u_{\sigma}\in U_{\sigma}\;$ and $\;u'_{\sigma'}\in U_{\sigma'}\;,$
$u_{\sigma}\sim u'_{\sigma'}$ if and only if $\phi_{M^t(M'^{-1})^t} $ is bi-regular on $u_{\sigma}$ and 
$\phi_{M^t.(M'^{-1})^t}(u_{\sigma})=u'_{\sigma'}$.

The quotient of $\mathfrak{C}$ under this equivalence relation is a smooth variety, $X_{\Sigma}$, called the \textbf{toric variety} associated to $\Sigma$. Each $U_\sigma$ is a \textbf{chart} of $X_{\Sigma}$.

For each $\sigma = Cone (M)\in \mathfrak{M}$ the mapping
\[
\begin{array}{cccc}
	\pi|_{U_\sigma}:
		& U_\sigma
			&\longrightarrow 
				&\K^n\\
		& x	&\mapsto
				& \phi_{M^t}(x)
\end{array}
\]
is a regular morphism, compatible with the gluing. 

The induced  regular morphism
\[
\pi:X_{\Sigma} \longrightarrow \K^n
\]
is called the \textbf{toric modification } associated to $\Sigma$. The morphism $\pi$ is a proper bi-rational morphism and it is bi-regular in the torus ${(\K^*)}^n$.

    Given a variety $V\subset \K^n$, the inverse image $\pi^{-1}(V)$ is called  the \textbf{total transform} of $V$ under the projection $\pi$.

    Let $V^*\subset V$ be the set of points in $V$ lying outside the coordinate hyperplanes, that is,
\[
V^{*}:=V-V(x_1\cdots x_n).
\] 
The \textbf{strict transform} of $V$ under the projection $\pi$ is defined by 
    \begin{eqnarray*}\label{transestr}
        \widetilde{V}:= \overline{\pi^{-1}(V^{*})},
    \end{eqnarray*} 
where $\overline{A}$ denotes the Zariski closure of $A$.

Given $\sigma\in \mathfrak{M}$, there is a natural action 
\[
\begin{array}{cccc}
        (\K^{*})^n \times U_{\sigma} &\longrightarrow & U_{\sigma}\\
         (\lambda_1,\ldots ,\lambda_n)\times (z_1,\ldots,z_n) & \mapsto& (\lambda_1 z_1,\ldots ,\lambda_n z_n)
    \end{array}
\]
that is compatible with the gluing in $X_\Sigma$.

Given $\lambda\subset\{ 1 ,\ldots ,n\}$ and $\sigma=Cone(M)\in\mathfrak{M}$, let $\tau$ be the cone generated by the $i^{th}$ columns of $M$ with $i\in\lambda$.

The {\bf orbit} associated to $\tau$ is the $(n-\#\lambda)$-dimensional torus given by
\[
{\mathcal O} (\tau ) := \{ (x_1,\ldots ,x_n) ; x_i=0 \text{ for } i\in\lambda \text{ and }  x_i\neq 0 \text{ for } i\notin\lambda\}.
\]
The set                                                                                                                      
${\mathcal O} (\tau )$ is well defined (it does not depend on the choice of $\sigma$) and is an orbit of the natural action on $X_\Sigma$.

\begin{Remark}\label{restringir a la orbita es como restringir al cono II}
Let $\sigma= Cone(M)\in\Sigma$ be of maximal dimension. Let $\tau$ be the cone in $\Sigma$ generated by the first $s$ columns of $M$. And let $h : U_\sigma\cap X_\Sigma \longrightarrow \K$ be a regular function.

 By Remark \ref{restingir es como poner ceros}, we have
\[
	h|_{{\mathcal O}(\tau )\cap U_\sigma}= h|_{Cone (e^{(s+1)},\ldots ,e^{(n)})}.
\]
\end{Remark}

\section{Newton non degenerate hypersurfaces.}

A hypersurface singularity  $\underline{0}\in H = V (f)$ is said to be {\bf Newton non degenerate} when, for any face $F$ of $NP (f)$, the hypersurface $V(f|_F)$ has no singularities outside the coordinate hyperplanes. In this section we will prove that a Newton non degenerate hypersurface singularity has a toric resolution. This result is proved in several texts (see for example \cite{Teissier:2004,Merle:1980}).

Let $\Sigma$ be a regular refinement of $\Sigma (f)$, let $\{\underline{0}\}\neq \tau$ be a cone in $\Sigma$, and let  $y$ be a point in the strict transform of $H$ under the toric modification given by $\Sigma$ with
\[
y\in \widetilde{H}\cap {\mathcal O} (\tau )\subset X_\Sigma.
\] 
Choose $\sigma = Cone (M)\in\mathfrak{M}$ such that the first $s$ columns of $M$ generate $\tau$.  We have
\[
y= (0,\ldots ,0,y_{s+1},\ldots , y_n)\in U_\sigma
\]
with $y_i\neq 0$ for $i\geq s+1$.

Since $Cone (M)$ is good for $f$, by Proposition \ref{El poliedro se convierte en el primer cuadrante trasladado}, $NP (f\circ\phi_{M^t})$ has only one vertex $\vertice$. By Remark \ref{Cuando el poliedro tiene solo unvertice es un monomio veces una unidad}, we have
\[
f\circ\phi_{M^t}= x^\vertice h
\]
with $h(\underline{0})\neq 0$.

Then 
\[
\widetilde{H}\cap U_\sigma = V(h)
\]
and, since $\underline{0}\notin V (h)$, $\tau\neq\sigma$.

By Remark \ref{restringir a la orbita es como restringir al cono II}
\begin{equation}\label{restringir a la orbita es como restringir al cono}
	h|_{{\mathcal O}(\tau )\cap U_\sigma }= h|_{Cone (e^{(s+1)},\ldots ,e^{(n)})},
\end{equation}
thus, the function $h|_{Cone (e^{(s+1)},\ldots ,e^{(n)})}$ does not depend on the first $s$ variables, then
\[
h (y)= h|_{Cone (e^{(s+1)},\ldots ,e^{(n)})}(1,\ldots ,1,y_{s+1},\ldots  ,y_n).
\]

By Proposition \ref{Como cambian las caras a las que hay que restringir},
\[
x^\vertice h|_{Cone (e^{(s+1)},\ldots ,e^{(n)})}= f\circ\phi_{M^t}|_{\vertice +Cone (e^{(s+1)},\ldots ,e^{(n)})}= f|_{face_\tau (f)}\circ\phi_{M^t}
\]

The hypersurface,  $V(f|_{face_\tau (f)})$ has no singularities outside the coordinate hyperplanes, then, $V(h|_{Cone (e^{(s+1)},\ldots ,e^{(n)})})$ is non singular in $(1,\ldots ,1, y_{s+1},\ldots ,y_n)$. That is, there exists $j\in\{ 1,\ldots ,n\}$ such that 
\[
\frac{\partial h|_{Cone (e^{(s+1)},\ldots ,e^{(n)})}}{\partial x_j} (1,\ldots ,1,y_{s+1},\ldots ,y_n)\neq 0
\]
since $h|_{Cone (e^{(s+1)},\ldots ,e^{(n)})}$ does not depend on the first $s$ coordinates, we have $j\geq s+1$.

By equation (\ref{restringir a la orbita es como restringir al cono})
\[
\frac{\partial h}{\partial{x_j}}(y)=\frac{\partial h|_{Cone (e^{(s+1)},\ldots ,e^{(n)})}}{\partial x_j} (y)=\frac{\partial h|_{Cone (e^{(s+1)},\ldots ,e^{(n)})}}{\partial x_j}(1,\ldots ,1,y_{s+1},\ldots y_n)\neq 0
\]
that is, $y$ is not a singular point of $\widetilde{H}$ and is transversal to ${\mathcal O}(\tau)$.

\section{The Gr\"{o}bner fan.}

The Gr\"{o}bner fan is the extension of the concept of dual fan to non principal ideals. All this may be computed using the program GFAN \cite{Gfan}.

Given $\omega\in (\R_{\geq 0})^n$, the $\mathbf{\omega-}$\textbf{initial part} of $f$ is
\[
In_{\omega}\;f:= f|_{face_\omega (f)}.
\]
That is, given $f$ as in (\ref{Pol}),
\[
	In_{\omega}\;f= \sum_{\{\mu\in\varepsilon(f); \omega \cdot \mu=\nu_{\omega}\;f\}} a_{\mu}x^{\mu}.
\]

\begin{Remark} \label{3.1.2}
Given $f \in \K[x_1,\ldots ,x_n]$, if $f (\underline{0}) \neq 0$ then $\nu_{\omega}(u)=0$ for all $\omega \in (\R_{\geq 0})^n$. And $In_\omega f= f(\underline{0})$ for all $\omega \in (\R_{> 0})^n$.
\end{Remark}

Let $I\subset \K [x_1,\ldots , x_n]$ be an ideal. Given $\omega\in (\R_{\geq 0})^n$, the \textbf{$\omega-$initial ideal} of $I$ is the ideal generated by 
the $\omega-$initial part of every polynomial in $I$, that is,
\begin{eqnarray}
\mathfrak{In}_{\omega}I := \langle In_{\omega}f; f\in I \;\rangle.
\end{eqnarray}

The variety $V(\mathfrak{In}_{(1,\ldots ,1)}I)$ is classically called the \emph{tangent cone} of $V(I)$. 

In the light of remark \ref{relacion de equiv en principales}, there is a natural way to associate a polyhedral fan to an ideal.
In $\R^n$, we define an equivalence relation as follows:\\
    $$\omega\sim \omega^{'}\;\Leftrightarrow \;\;\mathfrak{In}_{\omega}I= \mathfrak{In}_{\omega^{'}}I\;.$$
Given $\omega \in \R^n$, the closure in $\R^n$ of its equivalence class, that is, $$C_{\omega}(I):=\overline{\{u\in \R^n; \mathfrak{In}_{u}I= \mathfrak{In}_{\omega}I\}},$$
 is a polyhedral cone and the collection $$\Sigma^{'}(I):=\{C_{\omega}(I); \omega \in \R^n\}$$ forms a fan (\cite{MoraRobbiano:1988}).

Since we are interested in a local study of a variety, we will use the fan formed by the intersection of this fan and the first orthant, 
that is, 
\[
\Sigma(I)= \{\Sigma^{'}(I) \bigcap E_J;\; J\subset \{1,2,...,n\}\}\subset (\R_{\geq 0})^n
\]
 where $E_J=Cone(\{e^{(j)};\; j\in J\})$ 
and $e^{(j)}$ is the $j-$th element of the standard basis for $\R^n$.  We will refer to $\Sigma(I)$ as the \textbf{Gr\"{o}bner fan} of the ideal $I$.
The cones in the Gr\"{o}bner fan of $I$ will be called \textbf{Gr\"{o}bner cones} of $I$.

Given $f \in \K[x_1,\ldots ,x_n],$ we have
\[
\Sigma(f) =\Sigma(\langle f\rangle).
\]

\begin{Remark}\label{relacion de equivalencia en no principales}
	Given $\omega ,\omega'\in Int_{rel} \langle e^{(i_1)},\ldots ,e^{(i_s)}\rangle$ we have that  the initial ideals $\mathfrak{In}_{\omega}I$ and  $\mathfrak{In}_{\omega^{'}}I$ are equal if and only if $C_\omega (I) = C_{\omega'} (I)$.
\end{Remark}

\begin{Remark} \label{Remark4.4}
	Let $\sigma$ be a cone in $\Sigma(I)$ not contained in the coordinate hyperplanes. The cone $\sigma$ is of maximal dimension if and only if $\mathfrak{In}_{\sigma}I$ is a monomial ideal.
\end{Remark}

\begin{defi}\label{defi3.3.}
 We will say that $\sigma$ is a \textbf{good cone} for $I$ when it is contained in a  Gr\"{o}bner cone of $I$.
\end{defi}

Let $I \subseteq \K[x_1, x_2,\ldots ,x_n]$ be an ideal and let $\sigma$ be a good cone for $I$.
 We define 
\[
\mathfrak{In}_{\sigma} I := \mathfrak{In}_{\omega} I\qquad\text{where}\quad\omega\in Int_{rel}(\sigma).
\] 

\begin{Remark} \label{Parte inicial con respecto a caras}
	Let $I \subseteq \K[x_1, x_2,\ldots ,x_n]$ be an ideal, let $\sigma$ be a good cone for $I$ and let $\tau$ be a face of $\sigma$. Then, $\tau$ is a 	good cone for $I$ and
	\[
	\mathfrak{In}_{\sigma} \mathfrak{In}_{\tau} I =\mathfrak{In}_{\sigma} I.
	\]
\end{Remark}

\section{Reduced Gr\"{o}bner basis}

Given a term order, $\prec$, we define the  \textbf{initial term}, $In_{\prec}f$, of a non-zero polynomial $f\in \K [x_1,\ldots ,x_n]$, as
 the unique minimal term with respect to $\prec$.

Let $I\subset \K [x_1,\ldots ,x_n]$ be an ideal. The initial ideal of  $I$ with respect to $\prec$ is the ideal generated by the 
initial terms of the polynomials in $I$, that is, 
\[
\mathfrak{In}_{\prec}I:=\langle In_{\prec}f; f\in I \rangle.
\]

The closure of the equivalence class: \begin{eqnarray}C_{\prec}(I):=\overline{\{u\in \R^n;
 \mathfrak{In}_{u}I= \mathfrak{In}_{\prec}I\}}\end{eqnarray}
is a Gr\"{o}bner cone of maximal dimension.

\begin{prop}\label{prop4.20}
 Let $I\subset \K [x_1,\ldots ,x_n]$ be an ideal and let $\prec$ be a term order. There exist a set $G_{\prec}(I)=\{g_1, \ldots , g_r\}\subset I$ (called \textbf{reduced Gr\"{o}bner basis}) such that 
\begin{itemize}
       \item[\textit{i)}]  $C_{\prec}(I)$ is a good cone for $g_i$, $i \in \{1,\ldots ,r\}$. 
       \item[\textit{ii)}] $\mathfrak{In}_{\upsilon}(I) = \langle In_{\upsilon} (g_1), \ldots , In_{\upsilon}(g_r)\rangle $, for all $\upsilon \in C_{\prec}(I) $.
   \end{itemize}
\end{prop}

Proof.
\begin{itemize}
 \item[\textit{i)}] Direct consequence of [Prop 2.6, \cite{FukudaJensenThomas:2007}].
 \item[\textit{ii)}] Direct consequence of [Coro 2.14, \cite{FukudaJensenThomas:2007}].
\end{itemize}

\medskip

Given $\omega\in (\R_{\geq 0})^n$, we define a total order $\prec_{\omega}$ in $\K [x_1,\ldots ,x_n]$ by,

\begin{eqnarray}\label{eqnarray5} 
x^\alpha \prec_{\omega} x^\beta \Leftrightarrow \alpha \cdot \omega < \omega \cdot \beta \; \, or\,\,
\alpha \cdot \omega = \omega \cdot \beta \; \; and \; \; \alpha \prec_{lex} \beta,
\end{eqnarray}
where $\prec_{lex}$ is the lexicographical order.

\begin{prop}\label{prop4.19}
Let $I\subset \K [x_1,\ldots ,x_n]$ be an ideal. Given  $\sigma\in\Sigma(I)$  with $dim(\sigma) = n$  and $ \omega \in Int_{rel}(\sigma)$ we have that $\sigma = C_{\prec_{\omega}}(I). $
\end{prop}

Proof. Take $ \omega \in Int_{rel}(\sigma)$ and let ${\prec}_{\omega}$ be as in (\ref{eqnarray5}), then, by definition, we have 
$\mathfrak{In}_{{C_{\prec_\omega}}} (I) = \mathfrak{In}_{{\prec}_{\omega}} (I)$ and  $\mathfrak{In}_{\omega} (I)=\mathfrak{In}_{\sigma} (I).$ 
 Also, $  \mathfrak{In}_{{\prec}_{\omega}} (I)=\mathfrak{In}_{{\prec}_{\omega}} \mathfrak{In}_{\omega} (I)$ (by [Lemma 2.13,\cite{FukudaJensenThomas:2007}] ) 
and $  \mathfrak{In}_{{\prec}_{\omega}} \mathfrak{In}_{\omega} (I) = \mathfrak{In}_{\omega} (I)$ 
(by Remark \ref{Remark4.4}).\\
All above equations imply $\mathfrak{In}_{{C_{\prec_\omega}}}(I) = \mathfrak{In}_{\sigma} (I).$

\medskip

\begin{prop}\label{prop4.21}
Let $I\subset \K [x_1,\ldots ,x_n]$ be an ideal and let $\sigma$ be a good cone for $I$ of maximal dimension. There exists a system of generators $G_\sigma =\{ g_1,\ldots ,g_r\}$ of $I$ such that $\sigma$ is a good cone for 
$g_i$ for $i \in \{1,\ldots ,r\}$ and 
\[
\mathfrak{In}_{\upsilon}(I)= \langle In_{\upsilon} g_1,\ldots ,In_{\upsilon} g_r\rangle
\]
 for all ${\upsilon} \in \sigma.$
\end{prop}
Proof. Let $\sigma' \in \Sigma(I) $ be such that $\sigma \subset\sigma'$ (exists by definition of good cone). Take $\omega \in Int_{rel}\sigma'$  
and consider the term order $ {\prec}_{\omega}$. By Proposition \ref{prop4.19} we have 
\begin{eqnarray}\label{eqnarray50} 
\sigma \subset \sigma' = C_{\prec_{\omega}}(I).
\end{eqnarray}
By  Prop \ref{prop4.20}, there exists a system of generators $g_1,\ldots ,g_r$ of $I$, such that $C_{\prec_{\omega}}(I)$ is good for each $g_i.$ 
The result follows from (\ref{eqnarray50}).

\section{Tropical Variety}

\begin{defi}
 Given an ideal $I \subset \K [x_1,\ldots ,x_n],$  {\bf the tropical variety associated to $I$} is 
\[
\mathbf{TV}(I):=\{\omega\in \R^n ; \mathfrak{In}_{\omega}(I)\,\, \text{contains no monomial} \}.
\]
\end{defi}

\begin{prop}\label{Prop5.16100}
 An ideal $I\subset \K[x_1,\ldots,x_n]$ contains no monomials if and only if $V(I)\cap (\K^*)^n\neq\emptyset. $
\end{prop}

Proof. $\Leftarrow]$ If $ax^{\alpha} \in I$ and $z\in V(I),$ then $a z_1^{\alpha_1}\cdots z_n^{\alpha_n}=0$ so there exist $i$ such that $z_i=0$.
Hence, $z \notin (\K^*)^n.$\\
$\Rightarrow]$ Suppose that $V(I)\cap (\K^*)^n=\emptyset$. Then $V(I)\subset \bigcup _{i=1}^n \{x_i=0\}.$ 
The monomial $x^{(1,\ldots,1)}$ satisfies $V(x^{(1,\ldots,1)})\supset V(I)$ and, by the Nullstellensatz,  there exists $k$ such that 
$x^{(k,\ldots,k)}\in I$, i.e. $I$ contains a monomial.

\medskip

Proposition \ref{Prop5.16100} implies
\begin{eqnarray}\label{eqnarray500} 
\mathbf{TV}(I)=\{\omega\in \R^n ; V(\mathfrak{In}_{\omega}(I)) \cap (\K^*)^n\neq \emptyset\}.
\end{eqnarray}

\begin{prop}\label{Prop5.1611}
 Let $I\subset \K[x_1,\ldots,x_n]$ be an ideal. Given $\omega\in Int_{rel}(\sigma)$ with $\sigma\in \Sigma(I)$ a cone. We have
\[
\omega\in\mathbf{TV}(I) \Longleftrightarrow \sigma \subset \mathbf{TV}(I). 
\]
\end{prop}

Proof. 
\begin{itemize}
\item[$\Leftarrow$]  Suppose that $\sigma$ is contained in $\mathbf{TV}(I)$ then  $\omega \in Int_{rel}(\sigma)\subset \sigma$,
 implies $\omega\in\mathbf{TV}(I)$.
\item[$\Rightarrow$] Let $\sigma\subset \Sigma (I)$ be a cone. Let $\omega \in Int_{rel}(\sigma) $ be such that $\omega\in \mathbf{TV}(I)$.
\begin{itemize}
 \item Take $\omega' \in Int_{rel}(\sigma)$, by Remark \ref{relacion de equivalencia en no principales}
$\omega' \in \mathbf{TV}(I)$.
\item
Take $\nu \in \sigma \backslash Int_{rel}(\sigma)$. There exist $\tau \in \Sigma(I)$ face of $\sigma$ such that  $\nu \in Int_{rel}(\tau )$. By duality we have
\[
\varepsilon (In_{\omega}(f))=\varepsilon (In_{\sigma}(f)) \subset \varepsilon (In_{\tau}(f))=\varepsilon (In_{\nu}(f)).
\]
Then, if $In_\omega f$ is not a monomial neither is $In_\nu f$.
\end{itemize}
\end{itemize}

\medskip

From the last proposition, it follows that the tropical variety is a union of cones in the Gr\"{o}bner fan.

\begin{teo}{\bf Bieri-Groves} \cite{BieriGroves:1984,Sturmfels:1996}\label{Teorema de Bieri-Groves}
Let $\K$ be algebraically closed field and let $I \subset \K [x_1, \ldots , x_n]$ be a monomial-free prime ideal. Then the tropical variety $\mathbf{TV}(I)$ can be written as a  finite union of polyhedra of dimension
$dim\,V(I)$.
\end{teo}

\begin{prop}\label{las dimensiones de las tropicalzaciones son iguales}
Let $\K$ be an algebraically closed field and let $V=V(I)\subset\K^n$ be a pure dimensional variety not contained in the coordinate hyperplanes.
Given $\omega \in \mathbf{TV}(I)$, we have 
\[
 dim\, \mathbf{TV}(I)= dim\, \mathbf{TV}(\mathfrak{In}_{\omega}I).
\]
\end{prop}

Proof.
Suppose that $\nu\notin \mathbf{TV}(I)$ then, there exists $f\in I$ such that $In_\nu f$ is a monomial and, then, $In_\nu In_\omega f$ is also a monomial. Since $In_\omega f$ is in $\mathfrak{In}_{\omega}I$ we have $\nu\notin \mathbf{TV}(\mathfrak{In}_{\omega}I)$. This gives the inclusion
\[
\mathbf{TV}(\mathfrak{In}_{\omega}I)\subset \mathbf{TV}(I).
\]
Now, let $d$ be the dimension of $\mathbf{TV}(I)$. Given $\omega$ in the tropical variety $\mathbf{TV}(I)$, by Theorem \ref{Teorema de Bieri-Groves}, there exists a cone 
$\sigma$ that is a good cone for $I$ with $\omega\in\sigma$ and $dim\, \sigma= d$. 
 By Remark \ref{Parte inicial con respecto a caras}, $\mathfrak{In}_{\sigma}\mathfrak{In}_{\omega}I=\mathfrak{In}_{\sigma}I$. Since $\mathfrak{In}_{\sigma}I$ is not a monomial we have $\sigma\subset\mathbf{TV}(\mathfrak{In}_{\omega}I)$.

Then
\[
dim\,\mathbf{TV}(\mathfrak{In}_{\omega}I)\geq \mathbf{TV}(I).
\]
And the result is prooved.

\medskip

\begin{prop}\label{Prop5.161} 
Let $\K$ be an algebraically closed field and let $V=V(I)\subset\K^n$ be a pure dimensional variety not contained in the coordinate hyperplanes.
Given $\omega \in \mathbf{TV}(I)$, we have 
\[
dim \,(V\cap {\K^*}^n) = dim\, \mathbf{TV}(I)= dim\, \mathbf{TV}(\mathfrak{In}_{\omega}I).
\]
\end{prop}

Proof.
Direct consequence of Theorem \ref{Teorema de Bieri-Groves} and Proposition \ref{las dimensiones de las tropicalzaciones son iguales}.

\medskip

\section{Newton non-degenerate varieties.}\label{Variedades no degeneradas}

In this section we will extend the concept of Newton non degenerate variety to non-principal ideals. For complete intersections the concept was extended by Khovanskii \cite{Khovanskii:1977b,Khovanskii:1977} in 1976. 

\begin{defi}\label{definicion para interseccion completa}
Let $V=V(f_1,\ldots ,f_{k})\subset \K^n$ be a variety of dimension $n-k$. The variety $V$ is \emph{Newton non-degenerate} if 
for any $\omega\in (\R_{\geq 0})^n$ the variety 
\[
V(In_{\omega}f_1, In_{\omega}f_2,\ldots ,In_{\omega}f_{k})
\] 
is of dimension $n-k$ and has no singularities in $(\K^{*})^n$.
\end{defi}

Let $I$ be the ideal generated by $f_1,\ldots ,f_k$, then, $V(f_1,\ldots ,f_k)=V(I)$. Given an other set of generators of $I$, say, $f_1',\ldots ,f_k'$. The ideals
\[
\langle In_{\omega}f_1, In_{\omega}f_2,\ldots ,In_{\omega}f_{k}\rangle\quad\text{and}\quad \langle In_{\omega}f_1', In_{\omega}f_2',\ldots ,In_{\omega}f_{k}'\rangle
\]
are not necessarily the same. Definition \ref{definicion para interseccion completa} depends strongly on the generators of the ideal $I$ chosen.

The definition we propose extends the definition above to non complete intersection singularities and  does not depend on the generators.
\begin{defi}\label{nuestra definicion}
A singularity $\underline{0}\in V (I) \subset \K^n$, where $I\subset \K [x_1,\ldots ,x_n]$ is an ideal,  is \textbf{Newton non-degenerate} if for every $\omega\in (\R_{\geq 0})^n$, the
 variety defined by $\mathfrak{In}_{\omega}I$ does not have singularities in $(\K^{*})^n$. 
\end{defi}

To check if an ideal is Newton non degenerate or not, it is enough to check the condition for one vector in the relative interior of each cone of the Gr\"{o}bner fan that is contained in the tropical variety associated to $I$. These computations may be done using the software GFAN \cite{Gfan}.

Let  $I = \langle f_1, f_2, f_3\rangle \subset \C [x, y, z,w]$ where 
\[
f_1= xy + xw - yw,\,\, f_2 = xz- w^2\,\,\text{and}\,\, f_3 = yz - yw - w^2.
\]
The zero set $V (I)$ is a surface that has a rational singularity of multiplicity three (see \cite{LeTosun:2004} or \cite{Tjurina:1968}). 

In \cite{Mirna:2009} the Groebner fan and the different initial ideals are computed using GFAN \cite{Gfan}. Then, their singularities are computed using SINGULAR \cite{Singular} concluding that it is a Newton non degenerate singularity in the sense of Definition \ref{nuestra definicion}.

Let $V=V(f_1,\ldots ,f_{k})\subset \K^n$ be Newton non degenerate variety in the sense of Definition \ref{definicion para interseccion completa}. Let $\Sigma$ be a regular refinement of the dual fan $\Sigma(f_i)$ defined by $f_i$ for $i=1,\ldots ,k$ and let $\pi:X_{\Sigma}\longrightarrow \K^n$ be the torical modification associated to $\Sigma$. Let $\widetilde{V}$ be the strict transform of $V$. It is well-known \cite{Oka:1997} that the restriction $\pi:\widetilde{V}\longrightarrow V$ is a good resolution of $V$. The following sections are devoted to show the corresponding result for Newton non degenerate varieties in the sense of our definition (Definition \ref{nuestra definicion}).

\section{Tropical variety and strict transform}

Let $I \subset \K  [x_1,\ldots ,x_n]$ be an ideal, let $\Sigma$ be a regular refinement of the Gr\"{o}bner fan $\Sigma (I)$, and let $\widetilde{V(I)}$ be the strict transform of $V(I)$ under the modification
\[
\pi : X_\Sigma\longrightarrow \K^n.
\]

Let $\sigma=Cone (M)\in\Sigma$ be a cone of maximal dimension and
let $G_\sigma$ be a system of generators as in Proposition \ref{prop4.21}. 

We have
\begin{equation}\label{Los transformados de los generadores generan}
{\phi_{M^t}}_* (I) = \left< \{ g\circ \phi_{M^t}\}_{g\in G_\sigma}\right> .
\end{equation}

Take $g\in G_\sigma$.
Since $Cone (M)$ is good for $g$, by Proposition \ref{El poliedro se convierte en el primer cuadrante trasladado} and Remark \ref{Cuando el poliedro tiene solo unvertice es un monomio veces una unidad}, we have
\begin{equation}\label{pongo la transformada como producto de monomio por unidad}
g\circ\phi_{M^t}= x^\vertice h
\end{equation}
with $h(\underline{0})\neq 0$.

Set
\begin{equation}\label{Los generadores H}
	H_{G_\sigma}:= \{ h\in\K [x_1,\ldots ,x_n] \,;\, g\circ \phi_{M^t}=x^\vertice h\,\,\text{for some}\, g\in G_\sigma\,\text{and}\, h(\underline{0})\neq 0\}.
\end{equation}

We have
\[
{\phi_{M^t}}_*(I)\subset \langle H_{G_\sigma}\rangle\quad \text{and}\quad V({\phi_{M^t}}_*(I))\cap {\K^*}^n= V (H_{G_\sigma})\cap {\K^*}^n
\]
then
\begin{equation}\label{contenciones entre os ceros de las h y las transformadas}
	U_\sigma\cap \widetilde{V (I)}\subseteq V (H_{G_\sigma})
	\stackrel{(\ref{Los transformados de los generadores generan})}{\subseteq} 
	U_\sigma\cap\pi^{-1} (V(I)).
\end{equation}

In the following section we will see that, when $I$ is Newton non degenerate, $V(H_{G_\sigma})$ is non singular and transversal to $\pi^{-1}(\underline{0})$. As a consequence, the first inclusion is actually an equality.
\footnote{Since $G_\sigma$ is a Gr\"{o}bner base, it is in particular a normalised standard base and the equality holds always (See for example \cite{Hironaka:1964} or, for the analytic case, \cite{ArocaHironakaVicente:1977})}

\begin{prop}\label{Si intersecta esta en la variedad tropical}
	Let $I \subset \K  [x_1,\ldots ,x_n]$ be an ideal, let $\Sigma$ be a regular refinement of the Gr\"{o}bner fan $\Sigma (I)$, and let $\widetilde{V(I)}$ be the strict transform of $V(I)$ under the modification given by $\Sigma$.

Given $\tau\in\Sigma$ we have
\[
\widetilde{V(I)}\cap {\mathcal O}(\tau )\neq\emptyset \Rightarrow \tau\subset\mathbf{TV}(I).
\] 
\end{prop}

Proof.
Given $\tau\in\Sigma$,
take $\sigma=Cone (M)\in\Sigma$, of maximal dimension, such that $\tau$ is generated by the first $s$ columns of $M$. 

 Take $y\in\widetilde{V(I)}\cap {\mathcal O}(\tau )$, we have
\[
y= (0,\ldots ,0,y_{s+1},\ldots ,y_n)\in U_\sigma,\quad \text{with}\quad y_i\neq 0\quad\text{for}\quad i\in\{ s+1,\ldots ,n\}.
\]
Set
\[
y'= (1,\ldots ,1,y_{s+1},\ldots ,y_n)\in U_\sigma\cap {(\K^*)}^n.
\]
Let $G_\sigma$ be a system of generators as in Proposition \ref{prop4.21}. Let $H_{G_\sigma}$ be as in (\ref{Los generadores H}). Take $g\in G_\sigma$ and $h\in H_{G_\sigma}$ with $g\circ\phi_{M^t}=x^\vertice h$. By Remark \ref{restingir es como poner ceros}, the function $h|_{Cone (e^{(s+1)},\ldots ,e^{(n)})}$ does not depend on the first $s$ variables, then
\[
h|_{Cone (e^{(s+1)},\ldots ,e^{(n)})}(y')=h|_{Cone (e^{(s+1)},\ldots ,e^{(n)})}(y) = h (y)\stackrel{(\ref{contenciones entre os ceros de las h y las transformadas})}{=}0.
\]
where the second inequality follows from Remark \ref{restringir a la orbita es como restringir al cono II}.


By Proposition \ref{Como cambian las caras a las que hay que restringir},
\[
x^\vertice h|_{Cone (e^{(s+1)},\ldots ,e^{(n)})}= g\circ\phi_{M^t}|_{\vertice +Cone (e^{(s+1)},\ldots ,e^{(n)})}= g|_{face_\tau (g)}\circ\phi_{M^t}
\]
then
\[
\phi_{M^t} (y')\in V(In_\tau (g))\quad\text{for all}\quad g\in G_\sigma.
\]

By Proposition \ref{prop4.21}
\[
\mathfrak{In}_\tau (I) = \left< \{ In_\tau g\}_{g\in G_\sigma}\right>.
\]
then 
\[
\phi_{M^t} (y')\in V(\mathfrak{In}_\tau (I)).
\]
Since $y'\in {(\K^*)}^n$ then $\phi_{M^t} (y')\in {(\K^*)}^n$,
and the implication follows from (\ref{eqnarray500}).

\section{Toric resolution}

\begin{teo}\label{Teorema1}

    Let $I\subset \K[x_1,\ldots ,x_n]$ be an ideal and let $\Sigma$ be a regular refinement of the Gr\"{o}bner fan of $I$. 
If $\underline{0}\in V=V(I)$ is a Newton non-degenerate singularity, then the strict transform $\widetilde{V}$ of $V$ under the
 toric modification 
\[
\pi:X_{\Sigma}\longrightarrow \K^n
\] 
associated to $\Sigma$ is non-singular.
\end{teo}

Proof.
Let $\Sigma$ be a regular refinement of $\Sigma (I)$, let $\{\underline{0}\}\neq \tau$ be a cone in $\Sigma$, and let  $y$ be a point in the strict transform of $V(I)$ under the toric modification given by $\Sigma$ with
\[
\widetilde{V(I)}\cap {\mathcal O} (\tau )\subset X_\Sigma.
\] 
Choose $\sigma = Cone (M)$ of maximal dimension such that the first $s$ columns of $M$ generate $\tau$.  We have
\[
y= (0,\ldots ,0,y_{s+1},\ldots , y_n)\in U_\sigma
\]
with $y_i\neq 0$ for $i\geq s+1$.

Set
\[
y'= (1,\ldots ,1,y_{s+1},\ldots , y_n)\in U_\sigma\cap {(\K^*)}^n.
\]

Let $G_\sigma=\{ g_1,\ldots ,g_r\}$ be a system of generators as in Proposition \ref{prop4.21}. The cone $\sigma$ is good for each $g\in G_\sigma$ and
\begin{equation}\label{numero}
\mathfrak{In}_\tau (I) = \left< \{ In_\tau g\}_{g\in G_\sigma}\right>.
\end{equation}

The set $G_\sigma$ is a system of generators of $I$, then 
\begin{equation}\label{Los transformados de los generadores generan}
{\phi_{M^t}}_* (I)= \left< \{ g\circ \phi_{M^t}\}_{g\in G_\sigma}\right> 
\end{equation}

Let $H_{G_\sigma}=\{h_1,\ldots ,h_r\}$ be as in (\ref{Los generadores H}). We have
\[
\widetilde{V}\cap U_\sigma \subseteq V(h_1,\ldots ,h_r).
\]

By Remark \ref{restringir a la orbita es como restringir al cono II},
\begin{equation}\label{restringir a la orbita es como restringir al cono III}
	h_i|_{{\mathcal O}(\tau )\cap U_\sigma}= h_i|_{Cone (e^{(s+1)},\ldots ,e^{(n)})}.
\end{equation}
The functions $h_i|_{Cone (e^{(s+1)},\ldots ,e^{(n)})}$ do not depend on the first $s$ variables, hence
\[
h_i (y)= h_i|_{Cone (e^{(s+1)},\ldots ,e^{(n)})}(y').
\]

By Proposition \ref{Como cambian las caras a las que hay que restringir},
\[
x^{\vertice_i}h_i|_{Cone (e^{(s+1)},\ldots ,e^{(n)})}= g_i\circ\phi_{M^t}|_{\vertice_i+Cone (e^{(s+1)},\ldots ,e^{(n)})}=In_\tau g_i\circ\phi_{M^t}.
\]

Continuing as in the proof of Proposition \ref{Si intersecta esta en la variedad tropical} $\tau\subset\mathbf{TV}(I)$. Then, by Proposition \ref{Prop5.161} 
\[
\begin{array}{ll}
	dim \,V 
		& = dim\,{\bf TV}(V) = dim\, {\bf TV}(\mathfrak{In}_\tau (I))\\
		& \displaystyle \stackrel{(\ref{numero})}{=}dim\, V  (In_\tau g_1,\ldots ,In_\tau g_r)\\
		& = dim\, V(h_1|_{Cone (e^{(s+1)},\ldots ,e^{(n)})},\ldots ,h_r|_{Cone (e^{(s+1)},\ldots ,e^{(n)})}).
\end{array}
\]

The variety $V  (In_\tau g_1,\ldots ,In_\tau g_r)$   has no singularities outside the coordinate hyperplanes, then, $V(h_1|_{Cone (e^{(s+1)},\ldots ,e^{(n)})},\ldots ,h_r|_{Cone (e^{(s+1)},\ldots ,e^{(n)})})$ is non singular in $y'$. That is, 
\[
Rank\left(\frac{\partial  h_j|_{Cone (e^{(s+1)},\ldots ,e^{(n)})}}{\partial x_i}\right) (y')= dim\, V.
\]
Since $ h_j|_{Cone (e^{(s+1)},\ldots ,e^{(n)})}$ does not depend on the first $s$ variables
\[
Rank\left(\frac{\partial  h_j|_{Cone (e^{(s+1)},\ldots ,e^{(n)})}}{\partial x_i}\right)_{j=1,\ldots ,r,\, i=s+1,\ldots ,n} (y')= dim\, V
\]
Now
\[
\begin{array}{c}
\left(\frac{\partial  h_j|_{Cone (e^{(s+1)},\ldots ,e^{(n)})}}{\partial x_i}\right)_{j=1,\ldots ,r,\, i=s+1,\ldots ,n} (y')\\
 = \left(\frac{\partial  h_j|_{Cone (e^{(s+1)},\ldots ,e^{(n)})}}{\partial x_i}\right)_{j=1,\ldots ,r,\, i=s+1,\ldots ,n} (y)\\
=\left(\frac{\partial  h_j}{\partial x_i}\right)_{j=1,\ldots ,r,\, i=s+1,\ldots ,n} (y).
\end{array}
\]

Then $V(H_{G_\sigma})$ is non singular and transversal to $\pi^{-1}(\underline{0})$. Since 
\[
V(H_{G_\sigma})\cap {(\K^*)}^n= \pi^{-1}(V(I))\cap U_\sigma\cap {(\K^*)}^n,
\]
 we have that
\[
\widetilde{V(I)}= V(H_{G_\sigma})
\]
and the result is proved.

\begin{teo}
	Let $I \subset \K  [x_1,\ldots ,x_n]$ be an ideal, let $\Sigma$ be a regular refinement of the Gr\"{o}bner fan $\Sigma (I)$, and let $\widetilde{V(I)}$ be the strict transform of $V(I)$ under the modification given by $\Sigma$.

Given $\tau\in\Sigma$ we have
\[
\widetilde{V(I)}\cap {\mathcal O}(\tau )\neq\emptyset \Leftrightarrow \tau\subset \mathbf{TV}(I).
\] 
\end{teo}

Proof.
One implication is Proposition \ref{Si intersecta esta en la variedad tropical}.

Given $\tau\in\Sigma$,
take $\sigma=Cone (M)\in\Sigma$ of maximal dimension such that $\tau$ is generated by the first $s$ columns of $M$. 
 
Let $G_\sigma$ be a system of generators as in Proposition \ref{prop4.21}. Let $H_{G_\sigma}$ be as in (\ref{Los generadores H}).  Take $y\in {(\K^*)}^n$ such that $In_\tau g(y)=0$ for all $g\in G_\sigma$.

For each $g\in G_\sigma$ write, the same way we did in (\ref{pongo la transformada como producto de monomio por unidad}),
\[
g\circ\phi_{M^t} = x^{\vertice_g}h_g\quad\text{with}\quad h_g (\underline{0})\neq 0.
\]
We have
\[
\widetilde{V(I)}=V(\{h_g\}_{g\in G_\sigma}).
\]

Set
\[
 z=(z_1,\ldots ,z_n)=\phi_{M^t}^{-1}(y),
\]
we have
\[
0= In_\tau g(z) =z^{\vertice_g}(h_g|_{Cone (e^{(s+1)},\ldots ,e^{(n)})})(z)
\]
then
\[
0= h_g|_{Cone (e^{(s+1)},\ldots ,e^{(n)})}(z) = h_g(0,\ldots ,0, z_{s+1},\ldots ,z_n).
\]
Therefore, we conclude
\[
(0,\ldots ,0, z_{s+1},\ldots ,z_n)\in {\mathcal O} (\tau)\cap \widetilde{V(I)}.
\]

\def\cprime{$'$} \def\cprime{$'$} \def\cprime{$'$}

\end{document}